\documentclass{article}
\usepackage{arxiv}

\usepackage[utf8]{inputenc}  
\usepackage[T1]{fontenc}  
\usepackage{hyperref} 
\usepackage{url}            
\usepackage{booktabs}       
\usepackage{amsfonts}       
\usepackage{nicefrac}       
\usepackage{microtype}      
\usepackage{lipsum}
\usepackage[blocks]{authblk}
\usepackage{amsmath,amssymb,amsfonts}
\usepackage{algorithmic}
\usepackage{graphicx}
\usepackage{multirow} 
\usepackage{textcomp}
\usepackage{tabulary}
\newcolumntype{K}[1]{>{\centering\arraybackslash}p{#1}}
\usepackage{indentfirst}
\usepackage{tikz}

\usetikzlibrary{positioning}
\usetikzlibrary{arrows.meta}
\usetikzlibrary{shapes}
\usetikzlibrary{arrows,shapes}
\usetikzlibrary{positioning}

\providecommand{\keywords}[1]{\textbf{\textit{Keywords---}} #1}

\title{Efficient Global Optimization using Deep Gaussian Processes}

\author{\textbf{Ali Hebbal} \footnote{Ph.D. student, ONERA, DTIS, Université Paris Saclay, Université de Lille, CNRS/CRIStAL, Inria Lille, ali.hebbal@onera.fr}}

\affil{ONERA, DTIS, Université Paris Saclay, Université de Lille, CNRS/CRIStAL, Inria Lille}
\author{\textbf{Loic Brevault}\footnote{Research Engineer, ONERA, DTIS, Université Paris Saclay, loic.brevault@onera.fr} and \textbf{Mathieu Balesdent} \footnote{Research Engineer, ONERA, DTIS, Université Paris Saclay, mathieu.balesdent@onera.fr}}
\affil{ONERA, DTIS, Université Paris Saclay, F-91123 Palaiseau Cedex, France}
\author{\textbf{El-Ghazali Talbi} \footnote{Professor at Polytech'Lille - University of Lille, el-ghazali.talbi@univ-lille1.fr} and \textbf{Nouredine Melab} \footnote{Professor at Polytech'Lille - University of Lille, nouredine.melab@univ-lille1.fr} }
\affil{ Université de Lille, CNRS/CRIStAL, Inria Lille, Villeneuve d'Ascq, France }

\begin{document}
\maketitle

\begin{abstract}
Efficient Global Optimization (EGO) is widely used for the optimization of computationally expensive black-box functions. It uses a surrogate modeling technique based on Gaussian Processes (Kriging). However, due to the use of a stationary covariance, Kriging is not well suited for approximating non stationary functions. This paper explores the integration of Deep Gaussian processes (DGP) in EGO framework to deal with the non-stationary issues and investigates the induced challenges and opportunities. Numerical experimentations are performed on analytical problems to highlight the different aspects of DGP and EGO.
\end{abstract}

\keywords{
Efficient Global Optimization, non-stationary Kriging, Deep Gaussian Processes, surrogate modeling.
}

\section{Introduction}
Bayesian algorithms are widely used to deal with expensive black-box function optimization. They are based on surrogate models, allowing the emulation of the statistical relationship between the design variables and the response (objective function and constraints), to predict its behaviour using a dataset also called Design of Experiments (DoE). The evaluation cost of the surrogate models is cheaper, so it is possible to evaluate a greater number of design candidates. A complete review on the surrogate models that are widely used in design optimization is given in \cite{wang2007review}. One of the most popular Bayesian optimization methods is "Efficient Global Optimization" (EGO) developed by Jones \textit{et al.} \cite{jones1998efficient}. It uses Kriging surrogate model \cite{stein2012interpolation} which is based on the Gaussian Process (GP) theory. The main advantage of Kriging is that in addition to the prediction, it provides uncertainty estimation of the surrogate model response. Based on these two outputs, infill criteria are constructed to iteratively add the most promising candidates to the dataset. These points are then evaluated on the expensive functions and the surrogate model is updated and so on, until a stopping criterion is satisfied.

Classical Kriging is a GP with a stationary covariance function, inducing a uniform smoothness of the prediction. While it is effective to approximate stationary functions, it causes a major issue in the prediction of non-stationary ones. Indeed, in many design optimization problems, the objective functions or constraints vary with a completely different smoothness along the input space, due to the abrupt change of a physical property for example. Different approaches have been proposed to overcome this issue. Direct formulation of non-stationary covariance function, is one of the most explored strategies. Higdon \textit{et al.} \cite{higdon1999non} introduced a non-stationary version of the squared exponential covariance function by convolving spatially-varying kernels. Paciorek \textit{et al.} \cite{paciorek2006spatial} extended this formulation to covariance functions that are positive definite in the Euclidean space and especially the Matern covariance function. However, these approaches are applicable only to a maximum of 3 dimensional problems \cite{paciorek2006spatial}. Instead of a direct formulation, another method is to use local stationary covariance functions to model non-stationary functions. For instance, Haas \cite{haas1990kriging} proposed a moving window approach where the training and prediction region move along the input space so the covariance function is considered stationary within this window, while Rasmussen and Gharmani \cite{rasmussen2002infinite} used different stationary GPs in different subspaces of the input space. These strategies present some limitations for the general case. First, in computationally expensive problems, data are sparsed and using a local surrogate model with sparser data may be problematic. Second, it may induce discontinuities at the boundaries of the subspaces. Another alternative strategy using non-linear mapping, has been introduced by Sampson and Guttorp \cite{sampson1992nonparametric} and consists in deforming the input space in order to model the non-stationary response by a stationary model. Following this work, Xiong \textit{et al.} \cite{xiong2007non}  proposed a sparse and flexible parametrization of the mapping function, by considering a piece-wise density function with parametrized knots, allowing to apply the non-linear mapping for high-dimensional problems.

This paper aims at investigating a new strategy to handle the non-stationary issue in EGO based on a promising surrogate modeling class that is Deep Gaussian Processes (DGP) \cite{damianou2013deep} and to raise the challenges and opportunities induced by the coupling of EGO and DGP.

The paper is organized as follows. First, a description of EGO is presented with a review on stationary GP and classical infill criteria used (Section~\ref{sec:2}). Then, DGPs are introduced, and their advantages over simple GPs are highlighted (Section~\ref{sec:3}). Next, the challenges that occur when combining DGP with EGO are discussed and a first Deep Efficient Global Optimization framework (DEGO) is presented (Section~\ref{sec:4}). Finally, a comparison on analytical test cases of the proposed approach (DEGO) with non-linear mapping developed by Xiong and stationary Kriging is performed (Section~\ref{sec:5}). 
\section{Efficient Global Optimization}
\label{sec:2}
\subsection{EGO framework}
\begin{figure}[t]
\centering
\centerline{\includegraphics[width=0.8\linewidth]{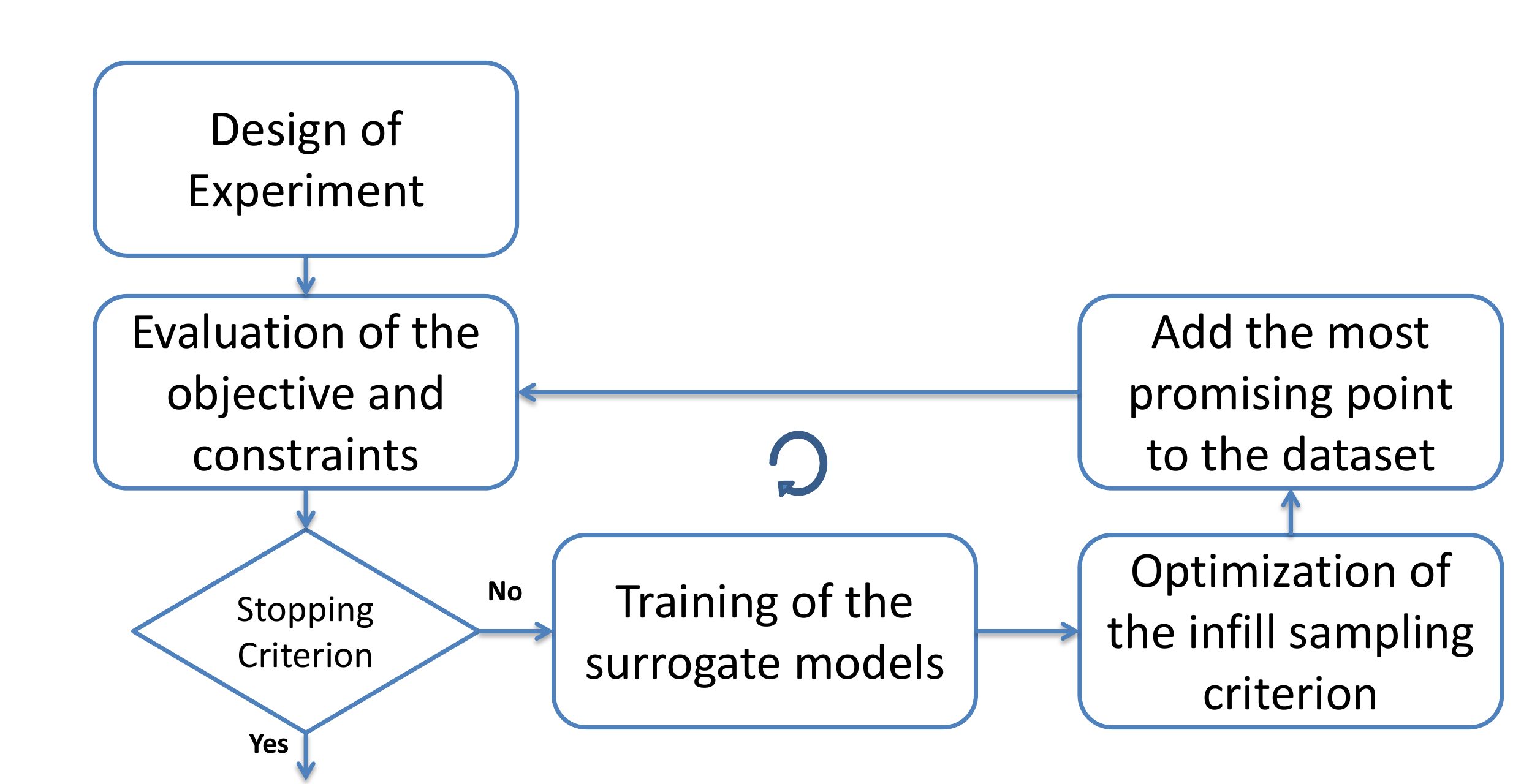}}
\caption{EGO framework}
\label{EGO}
\end{figure}

EGO is a Bayesian optimization algorithm dealing with expensive black-box optimization problems. It consists in sampling iteratively, using the prediction and uncertainty giving by the Kriging model, the most promising point based on an infill sampling criterion. This point is evaluated on the black box model and the surrogate-model is updated from the new training set, and a new point is sampled, and so on, until a stopping criterion is reached (Fig.~\ref{EGO}). Hence, the two important aspects in EGO are the Kriging surrogate model and the infill sampling criterion. 

\subsection{Kriging surrogate model}
A surrogate model is  built from a set of points called training set or design of experiment (DoE) $\mathcal{X}=\{\textbf{x}^{(1)},...,\textbf{x}^{(N)}\}$ and their associated response values $\mathcal{Y}=\{y^{(1)},...,y^{(N)}\}$ where $N$ is the number of observations and $\textbf{x}^{(i)} \in \mathbb{R}^d$. Then, it is possible to use this model to predict the output at a new point. The interest of a surrogate model is its cheap computational evaluation cost, thus instead of evaluating the expensive black-box function for a large number of points, it is possible to build a surrogate-model using fewer evaluations of the exact function, and predict its value in the other points using this model. 

The particularity of Kriging surrogate model is that it gives a variance estimation of the prediction, which makes it suitable in a global optimization framework.

A GP is used to describe a distribution over functions, it is a collection of infinite random variables, any finite number of which has a joint Gaussian distribution \cite{rasmussen2006gaussian}. It is defined by its mean function and covariance function. In GP regression, a GP prior is placed on the unobserved function (or latent) $f(\cdot)$ with a prior covariance function $k_\Theta(\textbf{x},\textbf{x'})$ that depends on a number of hyper-parameters $ \Theta$ and due to the fact that the trend of the response is \textit{a priori} unknown a constant mean function is considered (ordinary Kriging) \textit{i.e.}, $f(\textbf{x}) \sim \mathcal{N}(\mu,k_\Theta(\textbf{x},\textbf{x'}))$. Therefore, $f(\cdot)$ has a multivariate distribution on any finite subset of variables, in particular in $\mathcal{X}$ \textit{i.e.} $\textbf{f}|\mathcal{X} \sim \mathcal{N}(\textbf{1}\mu,\textbf{K}^\Theta_{NN})$ where $\textbf{K}^\Theta_{NN}$ is the covariance matrix constructed from the parametrized covariance function $k$ on $\mathcal{X}$ (further the dependence on $\Theta$ is dropped for notation simplicity). The choice of the covariance function is important, because it determines our prior assumptions of the function to be modeled.

A Gaussian noise variance is considered, such that the relationship between the latent function values $f(\mathcal{X})$ and the observed response $\mathcal{Y}$ is given by : $p(\textbf{y}|\textbf{f})= \mathcal{N}(\textbf{y}|\textbf{f},\sigma^2 \textbf{I})$. The marginal likelihood is obtained by integrating out the latent function $f(\cdot)$:

\begin{equation}
	p(\textbf{y}|\mathcal{X},\Theta) = \mathcal{N}(\textbf{y}|\mu, \textbf{K}_{NN}+\sigma^2 \textbf{I})
\end{equation}
Maximizing the marginal likelihood allows to train the GP by finding the optimal values of the hyper-parameters $\Theta, \mu$ and $\sigma$.
After the training, the prediction is made by considering a new point $\textbf{x}^*$ and using the conditional properties of a multivariate normal distribution \cite{rasmussen2006gaussian}:
\begin{equation}
	p(y^*|\textbf{x}^*,\mathcal{X},\mathcal{Y},\Theta)= \mathcal{N}(y^*|\hat{y^*},{\hat{s^*}}^2)
\end{equation}
with $\hat{y}^*$ the mean prediction and ${\hat{s^*}}^2$ the associated variance :
\begin{equation}
\hat{y^*}=\mu + \textbf{k}_{\textbf{x}^*}^T(\textbf{K}_{NN}^{-1}+\sigma^2 \textbf{I})^{-1} (\textbf{y}-\textbf{1}\mu)
\end{equation}
\begin{equation}
{\hat{s^*}}^2=k_{{\textbf{x}^*}{\textbf{x}^*}}-\textbf{k}_{\textbf{x}^*}^T(\textbf{K}_{NN}^{-1}+\sigma^2 \textbf{I})^{-1} \textbf{k}_{\textbf{x}^*}+\sigma^2
\end{equation}
where $k_{{\textbf{x}^*}{\textbf{x}^*}}=k({\textbf{x}^*},{\textbf{x}^*})$ and $\textbf{k}_{\textbf{x}^*}=[k(\textbf{x}^{(i)},{\textbf{x}^*})]_{i=1,...,N}. $
\subsection{Infill Sampling Criteria}
Different criteria have been developed for selecting infill sample candidates \cite{sasena2002flexibility}. They are based on a trade-off between exploration, by searching where prediction variance is high and exploitation by searching where prediction is minimized. The Probability of Improvement (PI) criterion samples the point where the probability of improving beyond the current minimum $y_{min}$ is maximum. That is $PI(\textbf{x})=P\left[y\leq y_{min}\right]=\Phi\left(\frac{y_{min}-\hat{y}(\textbf{x})}{\hat{s}(\textbf{x})}\right)$, where $\Phi(\cdot)$ is the Gaussian Cumulative Distribution Function (CDF). The higher values of $PI(\textbf{x})$ the higher chances that $\hat{y}(\textbf{x})$ is better than $y_{min}$. The inconvenient of this criterion is that only the probability is taken into account and not how much a point may improve the current best. This will add a lot of points around the current best point. 
To overcome the inconvenience of the Probability of Improvement, the Expected Improvement (EI) takes into account the improvement induced by a candidate that is defined as: $I(\textbf{x})=\text{max}\{0,y_{\text{min}}-f(\textbf{x})\}$. The EI is then computed as:
\par
\begin{eqnarray}
EI(\textbf{x})&=& \int{t . \text{PDF}_{I(\textbf{x})}(t)\text{d}t}\\
\label{eq:EI1} 
&=& (y_{\text{min}}-\hat{y}(\textbf{x}))\Phi\left(\frac{y_{\text{min}}-\hat{y}(\textbf{x})}{\hat{s}}\right)+\hat{s}\phi\left(\frac{y_{\text{min}}-\hat{y}(\textbf{x})}{\hat{s}}\right)
\label{eq:EI2} 
\end{eqnarray}%
\normalsize
where $\phi(\cdot)$ and $\Phi(\cdot)$ denote the Gaussian probability density function (PDF) and the Gaussian cumulative distribution function (CDF). The EI formula reveals two important terms. The first part is the same as in the PI, but multiplied by a factor that scales the EI value on the supposed improvement value. The second part expresses the uncertainty. It tends to be large when the uncertainty on the prediction is high. So, the EI is large for regions of improvement and also for regions of high uncertainty, allowing global refinement properties. \\
For constrained optimization problems, the feasibility of the constraints must be taken into account. Different techniques exist \cite{sasena2002flexibility}, optimization of the sub-problem (maximizing an unconstrained infill criterion) under approximated constraints (Direct method), under the constraints of an expected violation (Expected violation method), or maximizing the product of an unconstrained infill criterion with the probability of feasibility of the constraints (Probability of Feasibility method) are strategies widely used.
\section{Deep Gaussian Processes}
\label{sec:3}
\subsection{Introduction and Bayesian Training}
DGPs \cite{damianou2013deep} are a class of surrogate models based on the structure of neural networks, where each layer is a GP. It considers that the statistical relationship between the inputs and the response is expressed by a functional composition of GPs :
\begin{equation}
y=f_L(\textbf{f}_{L-1}(...(\textbf{f}_1(\textbf{f}_0(\textbf{x})))+\epsilon_L
\label{eq:}
\end{equation}
where $L$ is the number of layers and $\textbf{f}_l(.)$ is an intermediate GP. Each layer is composed of an input node $\textbf{h}_{l}$, an output node $\textbf{h}_{l+1}$ and a GP $\textbf{f}_l(.)$ mapping between the two nodes, getting the recursive equation: $\textbf{h}_{l+1}=\textbf{f}_l(\textbf{h}_{l})+\epsilon_l$ where $\epsilon_l$ is a Gaussian noise introduced in each layer (Fig.~\ref{DGP}). 

\tikzstyle{medium} = [circle, draw, thin, fill=blue!20, minimum height=2.5em]
\tikzstyle{textt} = [circle, thin, fill=white,align=center, minimum height=2.5em]
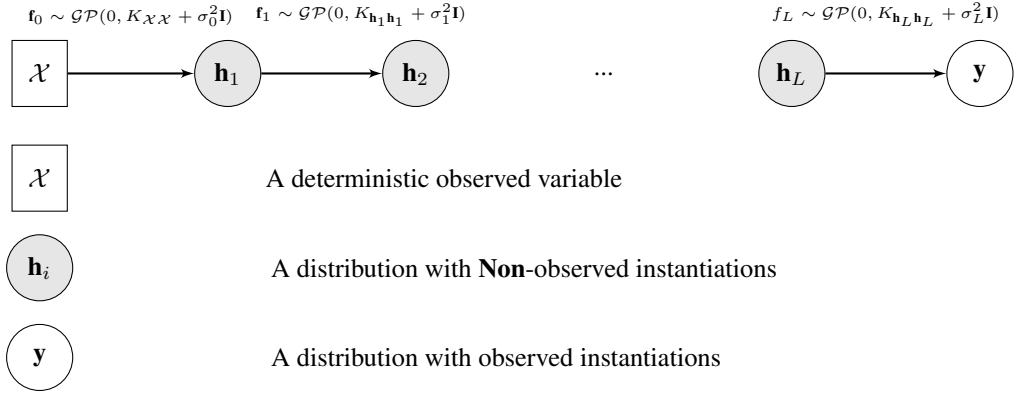
\begin{figure}[h]
\center
\begin{tikzpicture}[node distance=2.5cm, auto,>=latex', thick]
\path[use as bounding box] (2,-4) rectangle (10,1);
\node[medium,fill=white,rectangle,text width=0.5cm,align=center] (x) {${\mathcal{X}}$ };
\path[->] node[medium,fill=gray!20, right of= x] (h1) {$\textbf{h}_1$}
                  (x) edge node [above=0.5cm] {\tiny $\textbf{f}_0\sim \mathcal{GP}(0,K_{\mathcal{XX}}+\sigma_0^2 \textbf{I})$} (h1);
\path[->] node[medium,fill=gray!20, right of= h1] (h2) {$\textbf{h}_2$}
                  (h1) edge node [above=0.5cm] {\tiny \hspace{1cm} $\textbf{f}_1\sim \mathcal{GP}(0,K_{\textbf{h}_1 \textbf{h}_1}+\sigma_1^2 \textbf{I})$} (h2);                 
\path[->] node[textt,right of=  h2] (tt) {...};
            
\path[->] node[medium,fill=gray!20, right of= tt] (hN) {$\textbf{h}_L$};
\path[->] node[medium,fill=white, right of= hN] (y) {$\textbf{y}$}
                  (hN) edge node [above=0.5cm]  {\tiny $f_L\sim \mathcal{GP}(0,K_{\textbf{h}_L \textbf{h}_L}+\sigma_L^2 \textbf{I})$}  (y);
                  
\node[medium,fill=white, rectangle,text width=0.5cm,align=center,below=0.5cm of x] (x) {${\mathcal{X}}$ };
\node[ right =of x] (pp) {A deterministic observed variable};
\path[->] node[medium,fill=gray!20, below =0.3cm of x] (h1) {$\textbf{h}_i$};
\node[ right =of h1] (pp) {A distribution with \textbf{Non}-observed instantiations};
\path[->] node[medium,fill=white, below =0.3cm of h1] (y) {$\textbf{y}$};
\node[ right =of y] (pp) {A distribution with observed instantiations};
                  
\end{tikzpicture}
\caption{A representation of the structure of a DGP}
\label{DGP}
\end{figure}
The main difficulty induced by this cascade of GPs is that the intermediate nodes $\textbf{h}_l$ are latent variables \textit{i.e.} they are not observable as in a standard GP. Moreover, using non-linear covariance functions makes the overall composition no longer a GP. Hence, learning the hyperparameters of the model is challenging. To illustrate this, consider a DGP with one hiddden layer $\textbf{h}_1$. In this configuration a DGP is equivalent to a Gaussian Process Latent Variable Model (GPLVM)\cite{titsias2010bayesian}. The likelihood $p(\textbf{y}|\mathcal{X})$ is given by integrating over the latent variables $\textbf{h}_1$:

\begin{eqnarray}
p(\textbf{y}|\mathcal{X})&=&  \int p(\textbf{y}|\textbf{h}_1) p(\textbf{h}_1|\mathcal{X}) \text{d}\textbf{h}_1\\
&=& \int \mathcal{N}(0,K_{\textbf{h}_1 \textbf{h}_1}+\sigma_1^2\textbf{I})\mathcal{N}(0,K_{X X}+\sigma_0^2\textbf{I})\text{d}\textbf{h}_1\nonumber
\end{eqnarray}
for simplicity $\textbf{h}_1$ is taken one dimensional. The generalization to multi-dimensional hidden layer is done with the assumption of independence between dimensions.
Unfortunately, the integrals of Gaussians with respect to non-linear kernel functions are analytically intractable. To overcome this issue, the variational Bayesian approach is used. It consists of approximating the posterior distribution of the latent variables $p(\textbf{h}_1|\mathcal{Y})$ by a variational distribution $q(\textbf{h}_1)$ which has a factorized Gaussian form over the latent variables $\prod_{i=1}^N \mathcal{N}(\textbf{h}_1^{(i)}|\textbf{m}_1^{(i)},S_1^{(i)})$, where $\{\textbf{m}_1^{(i)},S_1^{(i)}\}_{i=1}^N$ are the variationnal parameters. Using Jensen's inequality, a first variational lower bound of the log likelihood is obtained: 
\par
\begin{eqnarray}
\log p(\textbf{y}|\mathcal{X}) &\geq & \int q(\textbf{h}_1) \log \frac{p(\textbf{y}|\textbf{h}_1) p(\textbf{h}_1|\mathcal{X})} {q(\textbf{h}_1)}\text{d}\textbf{h}_1\\
\label{eq1}
&=&   \int q(\textbf{h}_1) \log p(\textbf{y}|\textbf{h}_1)\text{d}\textbf{h}_1 - \int q(\textbf{h}_1) \log\frac{p(\textbf{h}_1|\mathcal{X})}{q(\textbf{h}_1)} \text{d}\textbf{h}_1  \nonumber \\
&=&  \tilde{F}(q)-KL(q||p) \nonumber
\end{eqnarray}
\normalsize
where the second term is the Kullback\text{-}Leibler (KL) divergence between the variational posterior distribution $q(\textbf{h}_1)$ and the prior over the latent variables $p(\textbf{h}_1)=\mathcal{N}(\textbf{h}_1|\textbf{0},K_{XX})$. Since both distributions are Gussian, it is computed analytically. However, the first term $ \tilde{F}(q)$  is still intractable, due to the fact that there is still the integration over non linear covariance function $K_{\textbf{h}_1 \textbf{h}_1}$. To deal with this problematic, Titsias and Lawrence \cite{titsias2010bayesian} introduced inducing variables that consists in augmenting with additional input-output pairs $\mathcal{Z}\in \mathbb{R}^{M\times d},\textbf{u} \in \mathbb{R}^{M}$, the latent space. Originally, inducing inputs were used as part of sparse GP to faster the GP regression $M<<N$. The direct methods of sparse GP \cite{snelson2006sparse} modify the GP prior in order to use the inverse of the covariance matrix of the inducing inputs $K_{MM}$ instead of $K_{NN}$. Titsias in \cite{titsias2009variational} used the inducing inputs in a variational framework to avoid changing the GP prior in the sparse approximation. Titsias and Lawrence \cite{titsias2010bayesian} used this same framework to overcome the intractability of $ \tilde{F}(q)$. This approach consisted of augmenting the joint distribution with a distribution over the inducing variables $p(\textbf{u})=\mathcal{N}(\textbf{u}|\textbf{0},K_{MM})$  to obtain :
\begin{eqnarray}
p(\textbf{y}|\textbf{h}_1) &=& \int_{\textbf{u},\textbf{f}_1} p(\textbf{y},\textbf{u},\textbf{f}_1|\textbf{h}_1) \text{d}\textbf{u} \text{d} \textbf{f}_1\\ &=& \int_{\textbf{u},\textbf{f}_1}p(\textbf{y}|\textbf{f}_1)p(\textbf{f}_1|\textbf{u},\textbf{h}_1)p(\textbf{u}) \text{d}\textbf{u} \text{d} \textbf{f}_1  \nonumber
\end{eqnarray}
The dependence on $\mathcal{Z}$ is dropped for notation simpicity. Then, the variational approach is used by approximating the true posterior of the latent variables $p(\textbf{f}_1,\textbf{u}|\textbf{y},\textbf{h}_1)=p(\textbf{f}_1|\textbf{u},\textbf{h}_1,\textbf{y})p(\textbf{u}|\textbf{y},\textbf{h}_1)$ with the variational distribution $q(\textbf{f}_1,\textbf{u})=p(\textbf{f}_1|\textbf{u},\textbf{h}_1)q(\textbf{u})$. The approximation $ p(\textbf{f}_1|\textbf{u},\textbf{h}_1,\textbf{y}) \approx p(\textbf{f}_1|\textbf{u},\textbf{h}_1) $  implies that $\textbf{u}$ is a sufficient statistic for $\textbf{f}_1$ which is true for an optimal set of $\mathcal{Z}$ and $\textbf{u}$, and $q(\textbf{u})$ is a free Gaussian distribution over the inducing variables. By using the same trick as in Eq.(9) a lower bound of log$[p(\textbf{y}|\textbf{h}_1)]$ that depends on  $\mathcal{Z}$, $q(\textbf{u})$ and the hyperparamters $\Theta_1$ of $\textbf{f}_1$ is obtained. Then, by optimizing this lower bound analyticly according to $q(\textbf{u})$ a tighter lower bound that depends only on $\mathcal{Z}$ and $\Theta_1$ is obtained. By replacing $\tilde{F}(q)$ in  Eq.(9) by this lower bound, an overall lower bound of the log likelihood log$[p(\textbf{y}|\mathcal{X})]$ that depends on $\mathcal{Z}$, $\Theta_0$, $\Theta_1$ and $q(\textbf{h}_1)$ is obtained. This lower bound is computable analyticly for kernels that are feasibly convoluted with the Gaussian density $q(\textbf{h}_1)$ such as the linear, the squared exponential and the Automatic Relevence Determination (ARD) squared exponential kernel.\\ The generalization to L number of layers consists in using the similar approximations in each layer with the assumption of independence between layers that is $q(\{\textbf{h}_l\}_{l=1}^L)=\prod_{l=1}^L q(\textbf{h}_l)$ \cite{damianou2013deep}. Hence, in a DGP configuration a lower bound of the likelihood is maximized according to $\{\Theta_l\}_{l=0}^L,\{\mathcal{Z}_l\}_{l=0}^L$ and $\{q(\textbf{h}_l)=\prod_{i=1}^N\mathcal{N}(\textbf{h}_l|\textbf{m}_l^{(i)},S_l^{(i)})\}_{l=0}^L$. 

\subsection{Other Training Approaches}
The main issue in the previous approach is the number of variational parameters $q(\textbf{h}_l)$ that increases linearly with the number of training datapoints, which complicates the optimization of the approximated likelihood. To overcome this issue, Dai \textit{et al.} \cite{dai2015variational} instead of considering the variational posteriors as individual parameters, take them as a transformation of observed data. Specificly, a recursive relationship links the variational parameters that is: $\textbf{m}_0 ^{(i)}=g_1(y^{(i)})$ and $\textbf{m}_l^{(i)}=g_l(\textbf{m}_{l-1}^{(i)})$ where $g_l$ is a multi-layer perceptron. This backprogration mechanism transforms the initialization of the variational parameters to the initialization of neural network parameters, which has been well studied in deep learning literature \cite{weymaere1994initialization}. Furthermore, the variational parameters are moved coherently during the optimization process. Bui \textit{et al.} \cite{bui2016deep} proposed a deterministic approximation for DGPs based on an approximated Expectation Propagation energy function, and a probabilistic back\text{-}propagation algorithm for learning. Salimbeni and Deisenroth \cite{salimbeni2017doubly} proposed a Doubly stochastic variational inference that does not assume the independence between layers and the form of the kernel functions, hence loosing the analytical tractability, that is bypassed through a crude Monte-Carlo sampling from the variational posterior.

\subsection{Prediction}
Once trained, the prediction using DGPs in a new point $\textbf{x}^*$ uses either normality assumptions in each layer or sampling strategies. The first approach consists of assuming that the predicted distribution at the layer $l-1$ is Gaussian. This Gaussian distribution is used as input in the layer $l$ and the mean and covariance of the non-Gaussian output distribution at the layer $l$ is obtained using the GPLVM prediction formula \cite{titsias2010bayesian}. This distribution is then assumed to be Gaussian to use the GPLVM prediction formula at layer $l+1$ and so on until reaching the final layer. The second approach uses sampling strategies (crude Monte-Carlo sampling) along the layers \textit{i.e.} using the mean and the covariance of the first layer, \textit{k} samples are generated following a Gaussian distribution, then each sample is passed to the next layer undergoing a new transformation following a Gaussian, and so on until reaching the last layer. So, $k$ values are obtained, the mean and standard deviation of these values are the predicted mean and its associated variance. The first method does not require sampling and then is faster than the sampling strategie approach, however it is based on assumptions which may impact the accuracy of the prediction.
\subsection{Advantage of DGPs over GPs}
\begin{figure}[h]%
\begin{minipage}[c]{0.45\linewidth}
\centering
\includegraphics[width=1\linewidth]{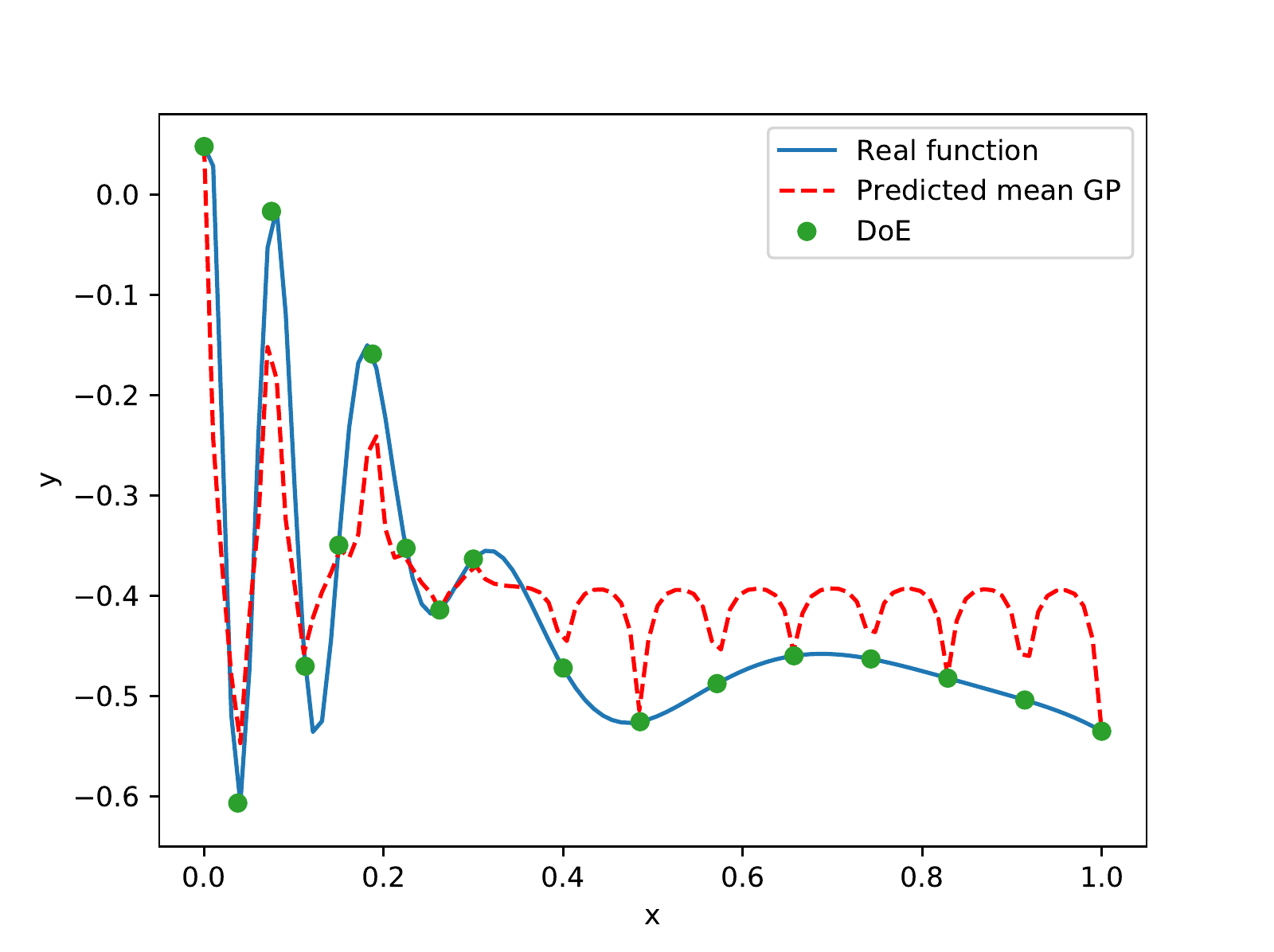}%
\caption{Regural GP}%
\label{regular}%
\end{minipage}
\hfill
\begin{minipage}[c]{0.45\linewidth}
\includegraphics[width=1\linewidth]{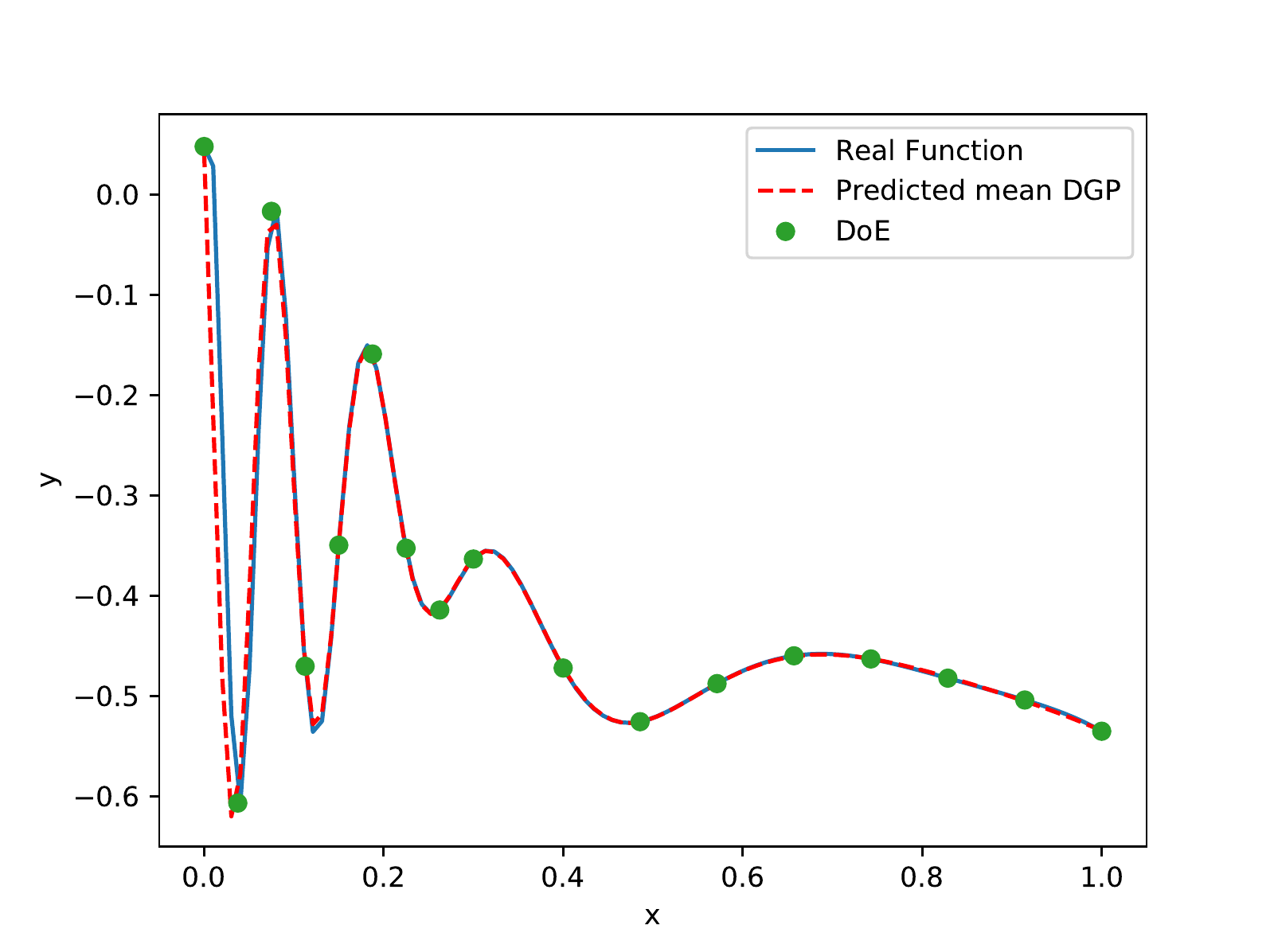}%
\caption{DGP}%
\label{dgpreg}%
\end{minipage}
\end{figure}

Since DGPs are a composition of GPs with different stationary kernels, the overall process is no longer a GP allowing the capture of non-stationarity (Fig. ~\ref{regular}, ~\ref{dgpreg}). Moreover it has been shown that DGPs handle scarce data and overcome the overfitting issue which is interesting in design optimization problem involving expensive black functions and induced uncertainty \cite{salimbeni2017doubly}.  

\section{DGPs and EGO}
\label{sec:4}
In this section the coupling of EGO and DGPs is discussed, highlighting the different arising challenges and opportunities. In \cite{dai2015variational}, an optimization experiment with DGPs was reported using the EI criterion. However, the combination EGO and DGP was not discussed, and as will be shown in this section the EI criterion cannot be used directly without prior assumptions.
\subsection{Number of induced inputs} 
In EGO, there is an iterative procedure based on an infill criterion to add a point to the dataset and update the surrogate model. However, in DGPs the number of induced inputs is directly depending on the cardinal of the dataset. So, in DEGO the number of induced inputs may change over the iterations. A small number of inducing inputs allows computational speed ups in the construction of the model by using reduced size matrices and less optimization parameters, while a high number allows accuracy in the approximation. Therefore, a trade-off between the two has to be set. Since in computationally intensive engineering design problems, the dataset is not large, the construction of the model is relatively cheap. One approach is to set the number of induced points in each layer equal to the number of data points over the iterations, another approach to advantage even more the accuracy is to set the number of induced variables in each layer equal to the sum of the number of data points and the number of infill points. However, since EGO is an iterative update of the surrogate model, it may be interesting to investigate robust approaches to vary the number of induced inputs considering the result of previous iterations. 
\subsection{Number of layers}
In DGPs, each layer transforms the output of the previous layer allowing the approximations of new features, and learning more complex representations. Experimentations with 2 to 5 layers were conducted in \cite{salimbeni2017doubly} and showed that the gain with increasing layers is achieved on very large dataset (one billion points). Since in engineering design problems, the dataset does not exceed the hundreds, two or three-layers configuration may be sufficient to catch the non stationarity of the model \cite{salimbeni2017doubly}. Nevertheless, an adaptive number of layers  and of the dimension of a layer, to EGO may be a promising strategy to explore. Indeed, one can begin with a two layer configuration in the first iterations, then, after a certain number of iterations or based on a given criterion that expresses the accuracy of the surrogate model, switches to a higher number of layers. 
\subsection{Training the model}
While in regular GP regression there is only kernel hyper-parameters to optimize in the training, in DGPs, in addition to kernel hyper-parameters in each layer, there are also the variational parameters. To faster the training, since in an EGO framework, the training is repeated in each iteration with an added point, it may be interesting after adding a consequent number of added points (to ensure a certain stability of hyper-parameters) to explore the use of the previous optimal configuration of the hyper-parameters as initialization for the next learning phase.
\subsection{Infill Criterion}
EGO is based on the balance between exploitation by searching where the value of the prediction is minimal, and exploration by searching where uncertainty is high. The EI is a widely used infill criterion. However, a direct application of the formula of the EI (Eq.~\ref{eq:EI2}) can not be computed analytically, since a DGP is not a Gaussian model. However, when the Gaussian assumption approach is used for prediction the direct formula of the EI can be used, since the approximated prediction is considered Gaussian.  On the other hand if the sampling approach is used for prediction, then, the expected improvement $EI(\textbf{x})$ is approximated using sampling on the value of the improvement $I(\textbf{x})$. In the constrained case, the probability of feasibility or the expected violation can also be approximated using a sampling strategy, or using the cumulative distributive function of a Gaussian with the assumption that the whole model behaves as a Gaussian model. So, The question that arises is when the assumption of Gaussian behavior is valid and if not how many samples must be used for an accurate approximation without useless computational cost ? 

\section{Experimentations}
\label{sec:5}
In this section, experimentations on two analytical test problems are performed to compare between EGO using DGPs (DEGO), standard EGO, and EGO using Xiong's non-linear mapping (NLEGO). The first problem is an unconstrained one dimensionnal optimization problem, and the second one is a constrained two dimensionnal problem (A Python implementation is publicly available \cite{DEGO2018}).

\begin{itemize}
\item In standard EGO, an ARD p-exponential kernel \cite{jones1998efficient} is used:  $k(\textbf{x},\textbf{x'})=l*\exp\{-\sum_{i=1}^D \theta_i (x_i -{x'}_i)^{p_i}\}$. The learning of the hyperparameters ($2$ hyperparameters by dimension) is done with CMA-ES \cite{hansen2006cma}.\\

	\item In NLEGO an ARD p-exponential kernel is used with the integration of the mapping $g(\cdot)$: $k(\textbf{x},\textbf{x'})=l*\exp\{-\sum_{i=1}^D \theta_i (g(x_i) -g({x'}_i))^{p_i}\}$. The learning of the hyper-parameters ($2+k_i$ by dimension where $k_i$ is the number of knots in dimension $i$ ) is done with CMA-ES.\\
	\item In DEGO, an ARD Gaussian kernel is used in each layer $k(\textbf{x},\textbf{x'})=l*\exp\{-\sum_{i=1}^D \theta_i (x_i -{x'}_i)^{2}\}$. A variational auto-encoded DGP with a $500$ dimensionnal encoder by layer is used for learning the hyperparameters  using l-BFGS-b optimization. The prediction is used with the Gaussian behaviour assumption. The learning and prediction of the DGP is performed using the toolbox PyDeepGP \cite{PyDeepGP2016}.
\end{itemize}

\subsection{1D unconstrained problem}
\subsubsection{Objective function}
The function to minimize is a one-dimensional non-stationary function presented in Equ.~\ref{eq:10}. It is a variant of the Xiong function \cite{xiong2007non} providing two-regions of interest in the minimization, one where the function varies with a high frequency $x\in [0,0.3]$ and the other where the function varies slowly $x \in [0.3,1]$ (Fig. ~\ref{regular}, ~\ref{dgpreg})
\begin{equation}
f(x)=-0.5(\sin[40(x-0.85)^4]\cos[2(x-0.95)]+0.5(x-0.9)+1), \text{  } x \in[0,1]
\label{eq:10}
\end{equation}
\subsubsection{Results}
\begin{figure}[h]%
\begin{minipage}[c]{0.45\linewidth}
\centering
\includegraphics[width=1\linewidth]{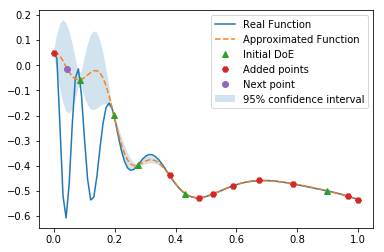}%
\caption{Best point maximizing the EI at iteration $k$ for DEGO-2}%
\label{iter1}%
\end{minipage}
\hfill
\begin{minipage}[c]{0.45\linewidth}
\includegraphics[width=1\linewidth]{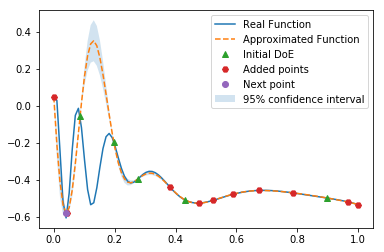}%
\caption{Best point added to the dataset at iteration $k+1$ for DEGO-2}%
\label{iter2}%
\end{minipage}
\end{figure}

20 initial DoE of five points are generated using a stochastic Latin Hypercube Sampling. 20 points are added using the EI criterion that is optimized with a differential evolution algorithm. \\
Table~\ref{table2} displays the mean best value attained ("DEGO l HL qD dynamic/m" corresponds to DEGO with l-hidden q-dimensional layers and a number of inducing inputs that is equal to the size of the dataset at each iteration if dynamic and equal to m otherwise) with the corresponding variance and the percentage of time observed to attain the true optimum that is $-0.60698$. Standard EGO is the less performing algorithm, which is expected since it does not take into account the non-stationarity of the function. NLEGO gives good results, since the function is divided in two different regions, four knots catch easily the non-stationarity. DEGO with one hidden layer gives similar results as NLEGO when the number of inducing inputs change over iterations, while setting the number of inducing points to 25 gives even better results. This is explained by the fact that more inducing inputs leads to a more accurate variational approximation. Moreover, the variance of the results on the 20 repetitions is lower for DEGO 1HL, which make it robust than the other algorithms. However, adding just one other hidden layer to the configuration gives poor results, this is due to the over-fitting induced by this complex configuration for a one dimensional function.\\
The interesting point in this comparison is that even for one dimensional problems where the non linear mapping is efficient, DEGO outperformed it with an adequate configuration. However, it also highlights the fact that an inadequate configuration has important impacts on the result.\\
Fig.~\ref{iter1},~\ref{iter2} illustrate respectively an iteration where the best point maximizes the EI and the next iteration where it is added to the dataset.

\begin{table}[h]
\centering	
\caption{Performance of the algorithms}

\begin{tabular}{|K{1.8cm}|K{1.8cm}|K{1.4cm}|K{1.4cm}|}
  \hline
  \textbf{Algorithm} & \textbf{The mean best value}&\textbf{ Variance} & \textbf{ \% of success} \\
  \hline
  EGO & -0.5493 & 0.00073 &  15\% \\
	\hline
	NLEGO 4 knots & \multirow{ 2}{*}{-0.5716} & \multirow{ 2}{*}{0.00115} & \multirow{ 2}{*}{45\%} \\
	\hline
	DEGO 1HL 2D dynamic  & \multirow{ 2}{*}{-0.5717} & \multirow{ 2}{*}{0.00127} & \multirow{ 2}{*}{50\%}\\
	\hline
	\textbf{DEGO 1HL 2D 25} & \multirow{ 2}{*}{\textbf{-0.581}} & \multirow{ 2}{*}{\textbf{ 0.00108 }} & \multirow{ 2}{*}{\textbf{55} \%}\\
	\hline
	DEGO 2HL 2D 25 & \multirow{ 2}{*}{-0.546} & \multirow{ 2}{*}{0.00067}  & \multirow{ 2}{*}{15 \%}\\
\hline

\end{tabular}

\label{table2}
\end{table}

\subsection{2D constrained problem}
\subsubsection{Objective function and constraint}

\begin{figure}[h]%
\centering
\includegraphics[width=0.7\linewidth]{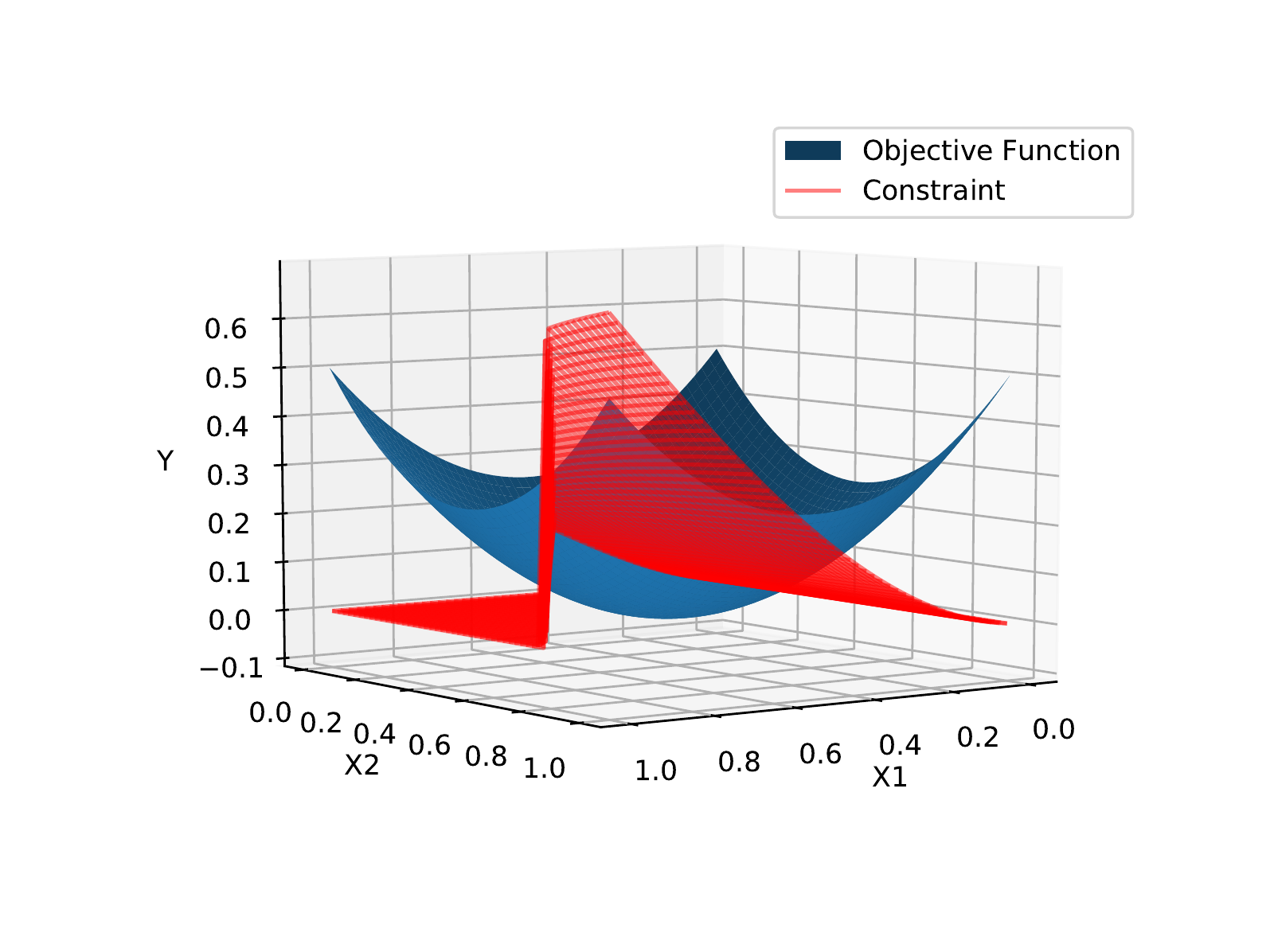}
\vspace{-0.8cm}
\caption{Objective and constraint functions 2D problem}
\label{exp2}
\end{figure}

The function to optimize is a simple two dimensional quadratic function: $f(x,y)=(x-0.5)^2+(y-0.5)^2$. While the constraint is non-stationary and feasible when equal to zero. An important discontinuity between the feasible and non feasible regions breaks the smoothness of the constraint (Fig.~\ref{exp2}). Therefore, the problem is challenging for standard GP, since the optimal region is exactly at the boundary of the discontinuity, requiring an accurate modelisation of the non-stationarity.

\subsubsection{Results}

20 initial DoE of 15 points are generated using a stochastic Latin Hypercube Sampling. 20 points are added using the EI criterion and the expected violation criterion for the constraint with a threshold of $10^{-3}$, and optimized with a differential evolution algorithm. Since the objective function is quadratic, a simple Gaussian process is used to approximate it in all the experimentations.

\begin{figure*}[!h]%
\centering
\includegraphics[width=1\linewidth]{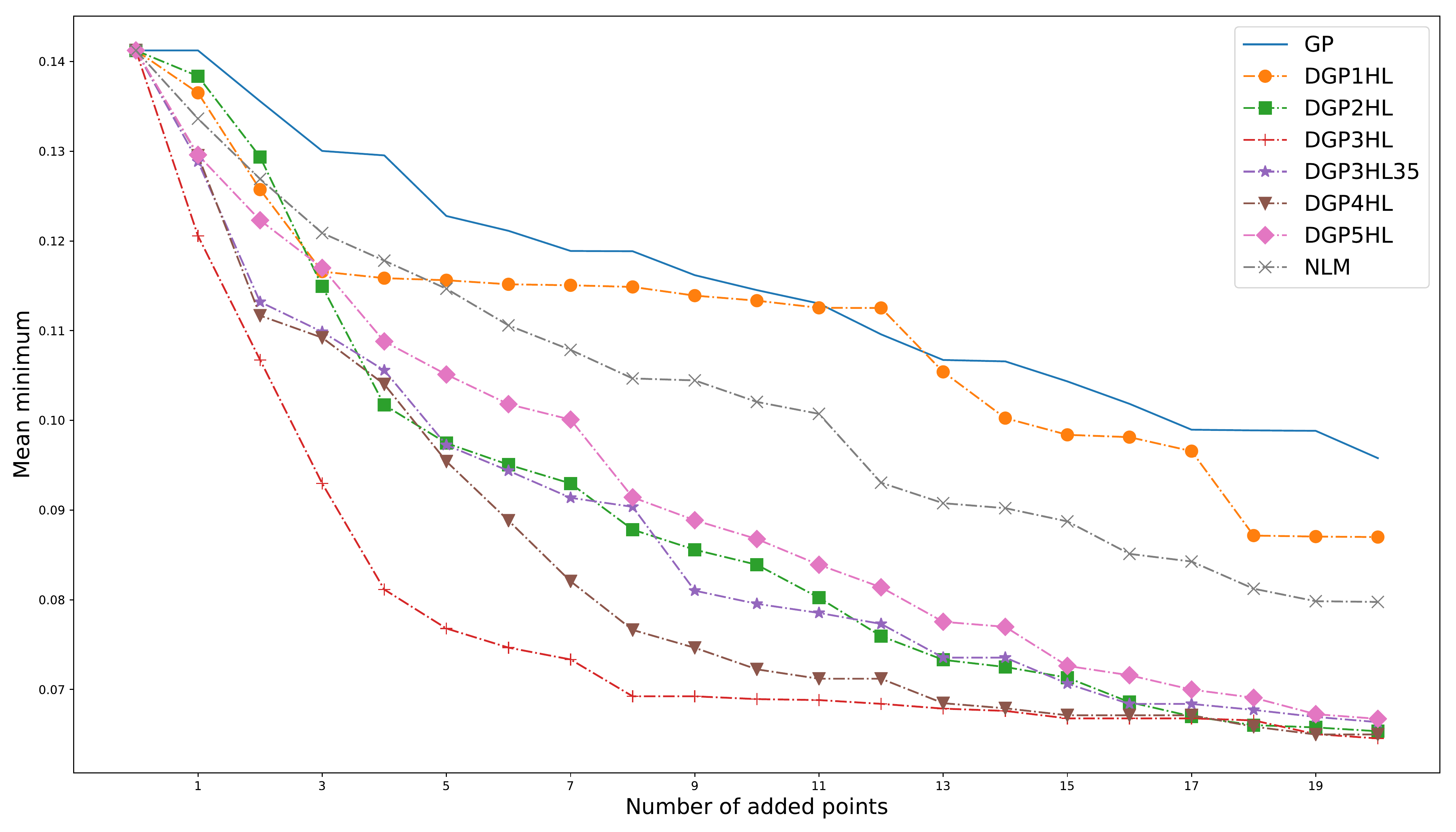}
\vspace{-0.5cm}

\caption{Evolution of the mean minimum according to the number of evaluations}
\vspace{-0.5cm}

\label{exp2evolution}
\end{figure*}

\begin{figure}[t]%
\begin{minipage}[c]{0.45\linewidth}
\centering
\includegraphics[width=1\linewidth]{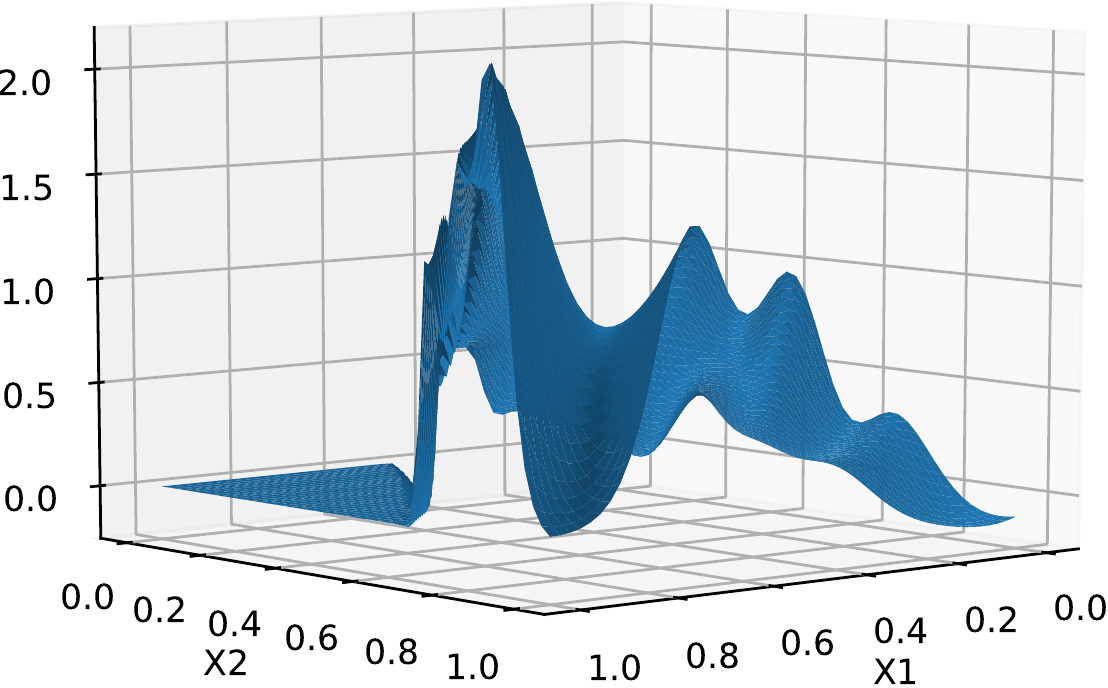}%
\caption{DGP approximation at the end of DEGO3 HL}%
\label{DGPcons}%
\end{minipage}
\hfill
\begin{minipage}[c]{0.45\linewidth}
\centering
\includegraphics[width=1\linewidth]{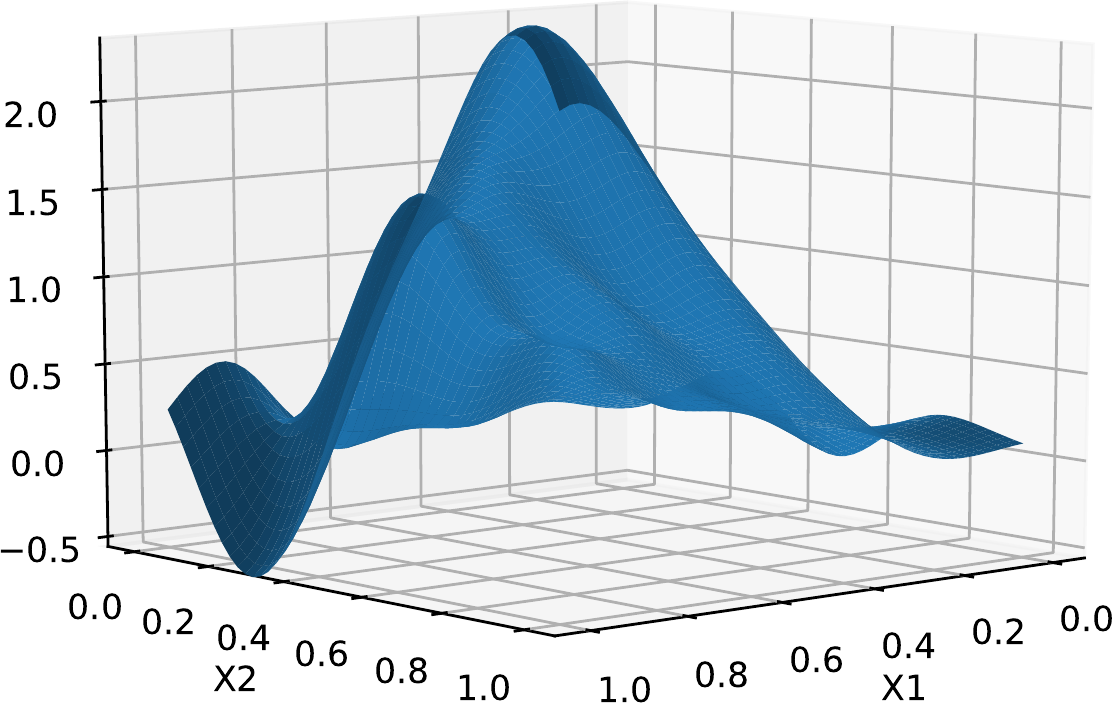}%
\caption{GP approximation at the end of EGO}%
\label{GPcons}%
\end{minipage}
\end{figure}

\begin{figure}[h]%
\centering
\vspace{-0.2cm}
\includegraphics[width=0.9\linewidth]{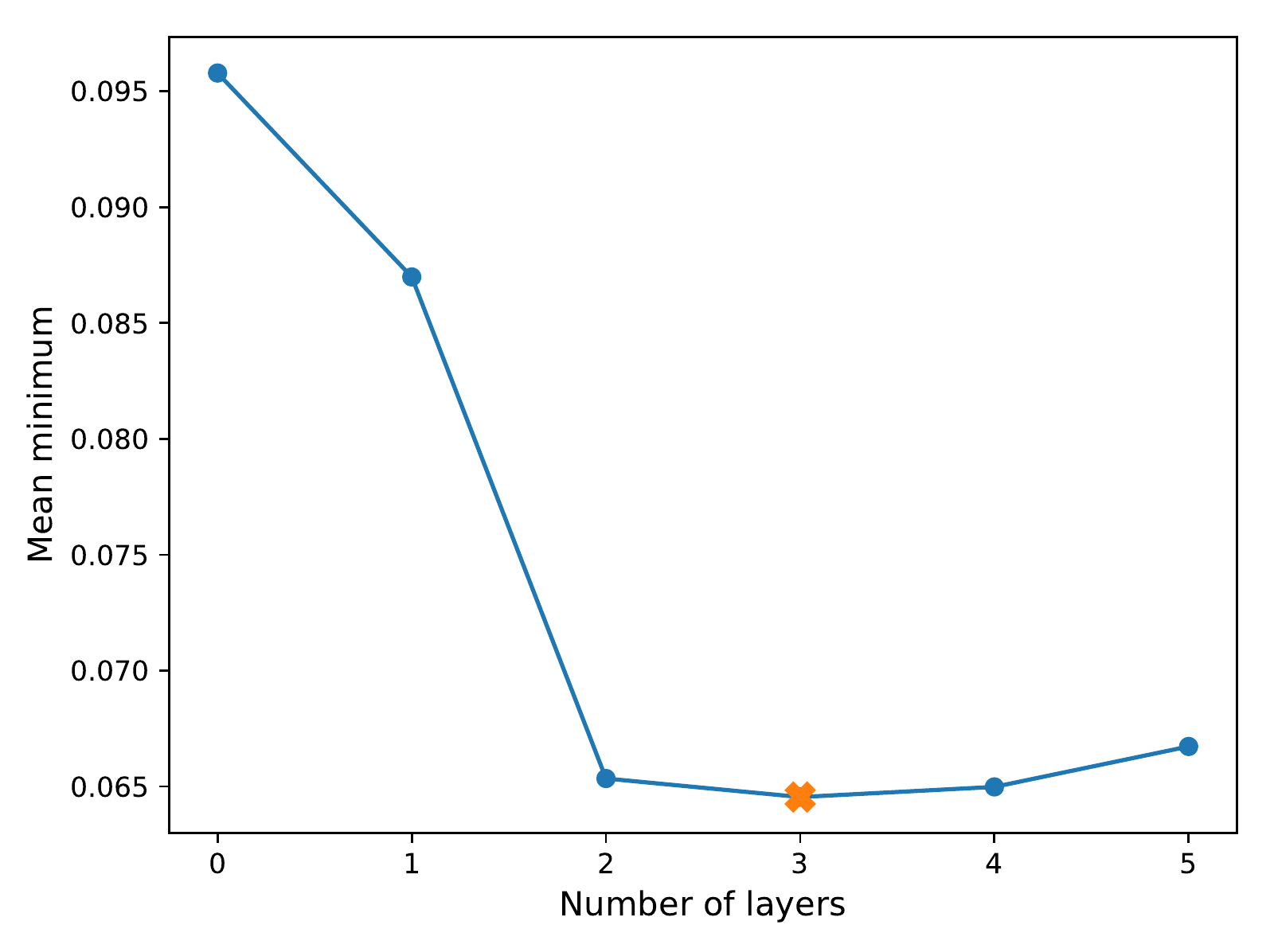}
\vspace{-0.5cm}
\caption{Evolution of the mean minimum according to the number of layers}
\vspace{-0.5cm}
\label{exp2layers}
\end{figure}

Table~\ref{table3} displays the mean best value attained by each algorithm with the corresponding variance and the percentage of time that it attains the true optimum that is $0.0602$. As in the first case standard EGO is not adapted to the problem, due to the discontinuity of the constraint. However, NLEGO does not performed as in the previous problem where its results were comparable to DEGO. This may be explained by two reasons. First, the non linear mapping is not well suited to the increase in the dimensionality of the problem, due to the fact that in the mapping the non-stationarity at a particular dimension may affect the other dimensions. Secondly, the non linear mapping can catch a change in the smoothness of the function in intervals but can not catch an abrupt change as the discontinuity in this constraint. DEGO with three layers and a dynamic number of induced points gives the better results. In contrast to standard GP, DEGO succeeds to capture the non stationarity of the model after adding a consequent number of points while classical EGO keeps a smooth modelisation and does not succeed to model accurately the discontinuity (Fig.~\ref{DGPcons},~\ref{GPcons}). The evolution of these algorithms according to the number of evaluations (Fig.~\ref{exp2evolution}) accentuates even more the superiority of DEGO and specially dynamic DEGO with 3 hidden layers. In fact, in the first iterations its speed of convergence is far more important than the other algorithms. Stopping the algorithms after adding 8 points would have given a more important gap between the different results. Another interesting aspect to observe is that the results given by DEGO are improved by adding hidden layers until reaching the three layers configuration. Adding more layers leads to a degradation of the results due to over-fitting (Fig.~\ref{exp2layers}). Finally, unlike the first problem, setting the number of induced points to a maximum does not give better results. All these numerical experiments illustrate the importance of the settings of the DGP configuration in EGO, but the appropriate configuration may provide better results than regular GP and non linear mapping.

\begin{table}[t]
\centering
\caption{Performance of the algorithms}
\begin{tabular}{|K{2cm}|K{2cm}|K{2cm}|K{2cm}|}
  \hline
  \textbf{Algorithm} & \textbf{The mean best value}&\textbf{ Variance} & \textbf{ \% of success} \\
  \hline
  EGO & 0.09579 & 0.001034 &  5 \% \\
  	\hline
	NLEGO 8 knots  & \multirow{ 2}{*}{0.07956}  & \multirow{ 2}{*}{0.000715} & \multirow{ 2}{*}{30 \%} \\
	\hline
	DEGO 1HL 10D dynamic  & \multirow{ 2}{*}{0.08699} & \multirow{ 2}{*}{0.000906} & \multirow{ 2}{*}{15 \%}\\
	\hline
	DEGO 2HL 10D dynamic & \multirow{ 2}{*}{0.065342} & \multirow{ 2}{*}{1.905 $10^{-5}$} & \multirow{ 2}{*}{50 \%}\\
	\hline
	\textbf{DEGO 3HL 10D dynamic }& \multirow{ 2}{*}{\textbf{0.06454 }}& \multirow{ 2}{*}{\textbf{1.32 $10^{-5}$}}  & \multirow{ 2}{*}{\textbf{75} \%}\\
	\hline
	DEGO 4HL 10D dynamic  &\multirow{ 2}{*}{0.06498} & \multirow{ 2}{*}{1.22 $10^{-5}$}  & \multirow{ 2}{*}{60 \%}\\
	\hline
	DEGO 3HL 10D 35  & \multirow{ 2}{*}{0.066358}  & \multirow{ 2}{*}{1.92 $10^{-5}$} & \multirow{ 2}{*}{45 \%}\\
	\hline
\end{tabular}

\label{table3}
\end{table}

\section{Conclusions}
\label{sec:6}
The coupling of EGO and DGP is a promising strategy to deal with optimization problems including non-stationary functions. In this paper, the challenges arising from the adaptation of EGO to DGP have been highlighted, and propositions were suggested in order to make the coupling possible. Finally, numerical experimentations confirmed the interest of DEGO giving the promising results, and also the difficulty in choosing the adequate configuration of the network. Hence, the necessity to provide an adaptive framework to set the configuration of the DGP according to the dynamic of EGO and the dimensions of the problem. Future works are to give a complete parametrization of DEGO (e.g. setting the number of layers, the width of the layer, the number of inducing points) according to an optimization problem, and to adapt a parallel infill criterion, allowing multiple points to be added in each iteration of DEGO, and to use this adaptive algorithm in real world computationally expensive black-box problems.
\bibliographystyle{unsrt}

\section{ACKNOWLEDGMENT}
This work was co-funded by ONERA-The French Aerospace Lab and Université de Lille/Inria Lille, in the context of a PhD thesis.

\bibliography{IEEEexample}

\begin{thebibliography}{10}

\bibitem{wang2007review}
Gary Wang and Songqing Shan.
\newblock Review of metamodeling techniques in support of engineering design
  optimization.
\newblock {\em Journal of Mechanical design}, 129(4):370--380, 2007.

\bibitem{jones1998efficient}
Donald~R Jones, Matthias Schonlau, and William~J Welch.
\newblock Efficient global optimization of expensive black-box functions.
\newblock {\em Journal of Global optimization}, 13(4):455--492, 1998.

\bibitem{stein2012interpolation}
Michael~L Stein.
\newblock {\em Interpolation of spatial data: some theory for kriging}.
\newblock Springer Science \& Business Media, 2012.

\bibitem{higdon1999non}
Dave Higdon, Jenise Swall, and J~Kern.
\newblock Non-stationary spatial modeling.
\newblock {\em Bayesian statistics}, 6(1):761--768, 1999.

\bibitem{paciorek2006spatial}
Christopher~J Paciorek and Mark~J Schervish.
\newblock Spatial modelling using a new class of nonstationary covariance
  functions.
\newblock {\em Environmetrics}, 17(5):483--506, 2006.

\bibitem{haas1990kriging}
Timothy~C Haas.
\newblock Kriging and automated variogram modeling within a moving window.
\newblock {\em Atmospheric Environment. Part A. General Topics},
  24(7):1759--1769, 1990.

\bibitem{rasmussen2002infinite}
Carl~E Rasmussen and Zoubin Ghahramani.
\newblock Infinite mixtures of gaussian process experts.
\newblock In {\em Advances in neural information processing systems}, pages
  881--888, 2002.

\bibitem{sampson1992nonparametric}
Paul~D Sampson and Peter Guttorp.
\newblock Nonparametric estimation of nonstationary spatial covariance
  structure.
\newblock {\em Journal of the American Statistical Association},
  87(417):108--119, 1992.

\bibitem{xiong2007non}
Ying Xiong, Wei Chen, Daniel Apley, and Xuru Ding.
\newblock A non-stationary covariance-based kriging method for metamodelling in
  engineering design.
\newblock {\em International Journal for Numerical Methods in Engineering},
  71(6):733--756, 2007.

\bibitem{damianou2013deep}
Andreas Damianou and Neil Lawrence.
\newblock Deep gaussian processes.
\newblock In {\em Artificial Intelligence and Statistics}, pages 207--215,
  2013.

\bibitem{rasmussen2006gaussian}
Carl Rasmussen and Christopher~KI Williams.
\newblock {\em Gaussian processes for machine learning}, volume~1.
\newblock MIT press Cambridge, 2006.

\bibitem{sasena2002flexibility}
Michael~J Sasena.
\newblock {\em Flexibility and efficiency enhancements for constrained global
  design optimization with kriging approximations}.
\newblock PhD thesis, University of Michigan Ann Arbor, MI, 2002.

\bibitem{titsias2010bayesian}
Michalis Titsias and Neil~D Lawrence.
\newblock Bayesian gaussian process latent variable model.
\newblock In {\em Proceedings of the Thirteenth International Conference on
  Artificial Intelligence and Statistics}, pages 844--851, 2010.

\bibitem{snelson2006sparse}
Edward Snelson and Zoubin Ghahramani.
\newblock Sparse gaussian processes using pseudo-inputs.
\newblock In {\em Advances in neural information processing systems}, pages
  1257--1264, 2006.

\bibitem{titsias2009variational}
Michalis Titsias.
\newblock Variational learning of inducing variables in sparse gaussian
  processes.
\newblock In {\em Artificial Intelligence and Statistics}, pages 567--574,
  2009.

\bibitem{dai2015variational}
Zhenwen Dai, Andreas Damianou, Javier Gonz{\'a}lez, and Neil Lawrence.
\newblock Variational auto-encoded deep gaussian processes.
\newblock {\em arXiv preprint arXiv:1511.06455}, 2015.

\bibitem{weymaere1994initialization}
Nico Weymaere and J-P Martens.
\newblock On the initialization and optimization of multilayer perceptrons.
\newblock {\em IEEE Transactions on Neural Networks}, 5(5):738--751, 1994.

\bibitem{bui2016deep}
Thang Bui, Daniel Hern{\'a}ndez-Lobato, Jose Hernandez-Lobato, Yingzhen Li, and
  Richard Turner.
\newblock Deep gaussian processes for regression using approximate expectation
  propagation.
\newblock In {\em International Conference on Machine Learning}, pages
  1472--1481, 2016.

\bibitem{salimbeni2017doubly}
Hugh Salimbeni and Marc Deisenroth.
\newblock Doubly stochastic variational inference for deep gaussian processes.
\newblock {\em arXiv preprint arXiv:1705.08933}, 2017.

\bibitem{DEGO2018}
{DEGO}.
\newblock {DEGO}: The python implementation of efficient global optimization
  with deep gaussian processes.
\newblock
  \url{https://github.com/M2CI-ONERA/M2CI-ONERA.github.io/tree/Deep-Gaussian-Process-EGO},
  since 2018.

\bibitem{hansen2006cma}
Nikolaus Hansen.
\newblock The cma evolution strategy: a comparing review.
\newblock In {\em Towards a new evolutionary computation}, pages 75--102.
  Springer, 2006.

\bibitem{PyDeepGP2016}
{PyDeepGP}.
\newblock {PyDeepGP}: The python implementation of deep gaussian processes.
\newblock \url{https://github.com/SheffieldML/PyDeepGP}, since 2016.

\end{thebibliography}

\end{document}